\documentclass[12pt]{article}
\usepackage{graphicx}
\usepackage[dvips]{epsfig}
\usepackage{subfigure}
\usepackage[a4paper]{geometry}
\usepackage{tikz-cd}
\usepackage{amsmath,amssymb,amsthm,amsfonts,hyperref,amscd}
\usepackage{times}
\usepackage{setspace}
\usepackage{verbatim}
\usepackage{cancel}
\usepackage{mathtools}
\usepackage[toc,page]{appendix}

\newtheorem{thm}{Theorem}[section]
\newtheorem{defn}[thm]{Definition}
\newtheorem{lem}[thm]{Lemma}
\newtheorem{prop}[thm]{Proposition}
\newtheorem{cor}[thm]{Corollary}
\newtheorem{rmk}[thm]{Remark}
\newtheorem{eg}[thm]{Example}
\theoremstyle{plain}
\theoremstyle{definition}

\counterwithin*{equation}{section}

\newcommand{\F}{\mathcal{F}}
\newcommand{\FF}{\mathbb{F}}
\newcommand{\W}{\mathbb{W}}
\newcommand{\x}{\chi_\mathcal{F}}
\newcommand{\dt}{\delta_T}
\newcommand{\db}{\delta_B}
\newcommand{\Db}{\Delta_B}

\newcommand{\kb}{\kappa_B}
\newcommand{\tb}{\tau_B}
\newcommand{\nD}{\nabla}
\newcommand{\V}{\text{vol}}

\newcommand{\addresses}{\bigskip\footnotesize

S. Hwang, \textsc{Department of Mathematics, Chung-Ang University, 84 HeukSeok-ro DongJak-gu, Seoul 06974, Republic of Korea.} \par\nopagebreak
\textit{E-mail address:}\texttt{seungsu@cau.ac.kr} 

\medskip

S. D. Jung, \textsc{Department of Mathematics, Jeju National University, Jeju 63243, Republic of Korea.} \par\nopagebreak
\textit{E-mail address:}\texttt{sdjung@jejunu.ac.kr} 

\medskip

J. Moon, \textsc{Department of Mathematics, Chung-Ang University, 84 HeukSeok-ro DongJak-gu, Seoul 06974, Republic of Korea.} \par\nopagebreak
\textit{E-mail address:}\texttt{dsfish999@cau.ac.kr} }

\title{Transverse Ricci solitons on a compact foliated manifold}
\author{Seungsu Hwang \and Seoung Dal Jung \and Jungwoo Moon}
\date{}

\begin{document}

\maketitle

\begin{abstract}
We investigate transverse Ricci solitons, the self-similar solutions of the transverse Ricci flow, on a compact foliated manifold. In particular, we show the relations between a taut Riemannian foliation and a transverse Ricci soliton. Moreover, we find some examples of transverse Ricci solitons.
\end{abstract}

\section{Introduction}

A Riemannian metric $g$ on a manifold $M$ is a \textit{Ricci soliton} if it is a self-similar solution of the Ricci flow. Equivalently, it is a metric $g$ on $M$ which satisfies the following: \begin{equation}Ric+\dfrac{1}{2}L_Xg=\lambda \text{ } g,\end{equation} where $\lambda$ is a constant, $Ric$ is the Ricci curvature, and $L_X$ is the Lie derivative with respect to a vector field $X$ on $M$. In particular, a \textit{gradient Ricci soliton} is a Ricci soliton satisfying $\displaystyle{X=\nD f}$ for some smooth potential function $f$ in (1.1). In general, a Ricci soliton may not be a gradient Ricci soliton. We refer to \cite{BD} for an example of Ricci soliton, which is not a gradient Ricci soliton. Nevertheless, we have the following well-known result on a compact manifold $M$ without boundary(shortly, a compact manifold $M$):
\begin{thm}\cite{Per}
Every Ricci soliton on a compact manifold is a gradient Ricci soliton.
\end{thm}

The idea of Theorem 1.1 depends on the nondecreasing property of Perelman's $\W$-functional, which is known as the following:
\begin{equation}
\W(f,g,\sigma)=\int_M [\sigma(S+|\nD f|^2)+f-n](4\pi\sigma)^{-\tfrac{n}{2}}e^{-f}dV,
\end{equation} where $f$ is a smooth function on a compact $M$, $S$ is the scalar curvature of the metric $g$, $\sigma$ is a time-dependent scale parameter, and $\displaystyle{n=\text{dim }M}$. To be specific, Theorem 1.1 has been proved since a gradient Ricci soliton is a critical point of (1.2).

On the other hand, there has been some research on the transverse analogy of the Ricci flow on a foliated manifold. Let $(M,\F,g)$ be a Riemannian manifold with a Riemannian foliation $\F$ and the normal bundle $Q$ of $\F$. Define a bundle-like metric $\displaystyle{g=g_{\F}\oplus g_Q}$ associated with $\F$ by $L_Ug_Q=0$ for any vector field $U$ tangent to $\F$, where $g_{\F}$ is the metric along $T\F$ and $g_Q$ is the metric on $Q$. Then the followings are known: First, Lovri\'c et al.\cite{LMR} introduced the transverse Ricci flow on $(M,\F,g)$ in 2000:
\begin{equation}
\begin{dcases}
\dfrac{dg_Q}{dt}=-2Ric^Q\\
\dfrac{dg_{\F}}{dt}=0,
\end{dcases}
\end{equation} where $Ric^Q$ is the transverse Ricci curvature of $(M,\F,g)$ with respect to $g_Q$.
Later, the analogy of Deturck's trick\cite{Det} of (1.3) was proved by Bedulli et al.\cite{BHV} and the $\FF^Q$ functional, the transverse analogy of Perelman's $\FF$ functional, is developed by Lin\cite{Lin}. Therefore, the self-similar solutions of (1.3) should be investigated for further research on the transverse Ricci flow.

Precisely, we show the following fundamental properties of the transverse analogy of Ricci solitons. Let a vector field $X$ on $(M,\F,g)$ be basic. That is, $[X,U]$ is tangent to $\F$ for any vector field $U$ tangent to $\F$. Then we have our first main result by the following:

\begin{thm}(cf. Theorem 3.6) Let $M$ be a compact Riemannian manifold with a Riemannian foliation $\F$ and let $g_Q(t)$ be a self-similar solution of (1.3) on $(M,\F)$. Then $g_Q(t)$ satisfies \begin{equation}
Ric^Q+\dfrac{1}{2}L_Xg_Q=\lambda \text{ } g_Q\end{equation} for all $t$, where $X$ is a basic vector field, $\lambda$ is a constant, and $Ric^Q$ is the transverse Ricci curvature on $(M,\F)$. Conversely, a family of transverse metric $g_Q(t)$ satisfying (1.4) is a self-similar solution of (1.3).
\end{thm}

The transverse metric $g_Q$ satisfying (1.4) is called a \textit{transverse Ricci soliton}. In particular, if $g_Q$ in Theorem 1.2 satisfies \begin{equation}
Ric^Q+\dfrac{1}{2}L_{\nD f}g_Q=\lambda \text{ } g_Q.\end{equation} for some basic function $f$ on $M$, then $g_Q$ is called a \textit{gradient transverse Ricci soliton}.

Inspired by Eminenti et al.\cite{ENM}, it is easy to obtain an analogous property of Theorem 1.1 within the context of transverse geometry on $(M,\F,g)$, where $M$ is compact and $\F$ is taut. i.e. there is a bundle-like metric $g$ which makes all leaves on $\F$ minimal. Therefore, it is natural to generalize Theorem 1.1 on an arbitrary Riemannian foliation on a compact manifold. 

Thus, we should introduce the following functional, which is the transverse analogy of $\W$ functional\cite{Per}, to show the aforementioned aim with using the same method of proof of Theorem 1.1:

\begin{equation}
\W^Q(f,g_Q,\sigma)=\int_M [\sigma(S^Q+|\nD f+\tb|^2)+f-q](4\pi\sigma)^{-\tfrac{q}{2}}e^{-f}dV,
\end{equation} where $f$ is a basic function on $(M,\F,g)$, $S^Q$ is the scalar curvature of $g_Q$, $\tb$ is the basic mean curvature of the foliation $\F$, and $\displaystyle{q=\text{dim } M- \text{dim }\F}$.

By the transverse Hodge-De Rham theorem\cite{Ton}, it is well-known that $\tb$ is a gradient of some basic function $h$ on $(M,\F,g)$ with compact $M$ if and only if $\F$ is taut. Thus, it is clear that $\W^Q$ is identified with $\W$ if a compact manifold $M$ is foliated by a taut Riemannian foliation $\F$ and hence (1.5) becomes a critical point of $\W^Q$. However, (1.5) is not a critical point of $\W^Q$ in general.  

Instead, the following is calculated as the critical point transverse metric of $\W^Q$(cf. Theorem 3.19):
\begin{equation}
Ric^Q+\dfrac{1}{2}L_{(\nD f+\tb)}g_Q=\lambda \text{ } g_Q,
\end{equation} Therefore, any transverse metric $g_Q$ satisfying (1.7) is also a transverse Ricci soliton, and such $g_Q$ is called a \textit{twisted gradient transverse Ricci soliton}. In particular, if $\tb$ is not a gradient vector field of a basic function $h$ on $(M,\F,g)$, then we call a twisted gradient transverse Ricci soliton \textit{nontrivial}. As mentioned in the above, (1.5) is identified with (1.7) if the given foliation $\F$ on a compact $M$ is taut. 

Hence, we have the following analogous result of Theorem 1.1 for the transverse Ricci solitons on $(M,\F,g)$ with compact $M$:
\begin{thm}
Let $M$ be a compact manifold of $(M,\F,g)$ and let $g_Q$ be a transverse Ricci soliton. i.e. $g_Q$ satisfies: $$Ric^Q+\dfrac{1}{2}L_Xg_Q=\lambda\text{ } g_Q$$ for some constant $\lambda$ and a basic vector field $X$. Then we obtain the following:
\begin{enumerate}
\item(cf. Theorem 4.10) If $\F$ is taut, then $g_Q$ is gradient.
\item(cf. Theorem 4.11) If $\F$ is non-taut, then $\lambda$ should be a negative constant. In other words, $g_Q$ is an expanding transverse Ricci soliton.
\item(cf. Corollary 4.13) Any transverse Ricci soliton cannot be a nontrivial twisted gradient transverse Ricci soliton.
\end{enumerate}
\end{thm}

\begin{rmk}\normalfont
The first variation formula of $\overline{\lambda^Q}$, which is defined by \begin{equation} \overline{\lambda^Q}=(\V(g))^{\tfrac{2}{q}}\text{inf}_{f \in \Omega_B^0(\F)}\int_M(S^Q+|\nD f+\tb|^2)e^{-f}dV,\end{equation} does not imply the nondecreasing property on a compact $(M,\F,g)$ along (1.3) unless $\F$ is taut\cite{Lin}(cf. Lemma 3.17). Thus, the authors do not know whether a transverse Ricci soliton on a compact $(M,\F,g)$, which is not gradient, exists or not. Nevertheless, we find an example satisfying the second statement of Theorem 1.3(cf. Example 5.2).
\end{rmk}

The content of this paper is structured as follows. In Section 2, we review the fundamental results of Riemannian foliations. In Section 3, we recall the previous works on the transverse Ricci flow to define transverse Ricci solitons, and Theorem 1.2 is proven. Moreover, we deal with the definition and properties of $\W^Q$ functional to introduce the (twisted) gradient transverse Ricci soliton. In Section 4, we prove the properties of (twisted) gradient transverse Ricci solitons on a compact foliated manifold to show Theorem 1.3. In Section 5, we give some examples of transverse Ricci solitons. In particular, we find an example of transverse Ricci solitons on a compact manifold with a non-taut Riemannian foliation.

\section{Preliminaries}

Let $M$ be a $(p+q)$-dimensional differential manifold with a foliation $\F$, let $T\F$ be the $p$-dimensional subbundle of the tangent bundle $TM$ tangent to $\F$. Define the transverse vector bundle $\displaystyle{Q:=TM/TF}$ of $\F$. That is, we consider the surjective $C^{\infty}$-linear map $\displaystyle{\pi: TM \xrightarrow[]{} Q}$ satisfying $\ker\pi=T\F$.

On a Riemannian manifold $M$, a foliation $\F$ is \textit{Riemannian} if the transverse part $g_Q:Q \times Q \xrightarrow{} \mathbb{R}$ of $g$ satisfies $\displaystyle{L_Ug_Q=0}$ for any vector field $\displaystyle{U \in \Gamma(T\F)}$, where $L_U$ is the Lie derivative with respect to $U$. In this case, the metric $g$ is bundle-like\cite{Rei}. We always denote $(M,\F,g)$ by a Riemannian manifold with a Riemannian foliation and the associated bundle-like metric. In particular, we always say \text{compact} $(M,\F,g)$ if $M$ is a compact manifold without boundary.

The transverse Levi-Civita connection on $(M,\F,g)$ is defined by\cite{Ton}

\begin{equation}
\begin{dcases}
\nD_E\pi(E')=\pi[E,\pi(E')] \text{ if }E \in \Gamma(T\F) \\
\nD_E\pi(E')=\pi(D_E\pi(E'))\text{ if }E \in \Gamma(Q)
\end{dcases}
\end{equation}
for any vector field $E'$ on $M$, where $D$ is the Levi-Civita conenction on $M$. We also use the notation $\nD_{tr}E'$ for the latter case of (2.1) if $E$ is undetermined. Let $\displaystyle{X, Y, Z, Z' \in \Gamma(Q)}$ be transverse vector fields. Then the \textit{transverse Riemann curvature tensor} $R^Q$ on $(M,\F,g)$ is defined by the following: \begin{equation}R^Q(X,Y)Z=\nD_Y\nD_XZ-\nD_X\nD_YZ+\nD_{[X,Y]}Z\end{equation} and the \textit{transverse Ricci curvature} $Ric^Q$(respectively, the \textit{transverse scalar curvature} $S^Q$) is defined by \begin{equation}
Ric^Q(X,Y)=\sum_{i=1}^qg_Q(R^Q(\mathbf{e}_i,X)\mathbf{e}_i,Y)
\end{equation} and
\begin{equation}
S^Q=\sum_{i=1}^qRic^Q(\mathbf{e}_i,\mathbf{e}_i),
\end{equation}
where $\{\mathbf{e}_i\}_{i=1,...,q}$ is a local transverse moving frame on $Q$ with respect to $g_Q$.
As the definition of transverse Riemann curvature and the definition of Riemann curvature tensor are parallel, the transverse second Bianchi identity is directly obtained by \begin{equation}
\nD_XR^Q(Y,Z)Z'+\nD_YR^Q(Z,X)Z'+\nD_ZR^Q(X,Y)Z'=0.
\end{equation}

From (2.5), we have
\begin{equation}
dS^Q=-2(\dt Ric^Q)
\end{equation} by taking the transverse trace with respect to $g_Q$ twice, where $\dt$ is the transverse divergence defined by
\begin{equation}
\dt=-\sum_{i=1}^qi_{\mathbf{e}_i}\nD_{\mathbf{e}_i}.
\end{equation}

We call $\F$ on $(M,\F,g)$ is an \textit{Einstein foliation} if its associated bundle-like metric $g$ satisfies
\begin{equation}
Ric^Q=\lambda g_Q
\end{equation} for some real constant $\lambda$. In particular, if $\lambda=0$, $\F$ is called a \textit{transversally Ricci-flat foliation}.

On the other hand, the mean curvature of the leaves of $\F$ is defined by \begin{equation}
\tau=\sum_{\alpha=1}^q \pi(D_{\xi_{\alpha}}\xi_{\alpha}),
\end{equation} where $\{\xi_{\alpha}\}_{\alpha=1,...,p}$ is the local orthonormal basis on $T\F$. Also, the mean curvature $1$-form given by $g_Q(\tau, \cdot)$ is denoted by $\kappa$. i.e. $\displaystyle{\tau=\kappa^{\sharp}}$. We call $\F$ \textit{minimal} if $\kappa$ vanishes and \textit{taut} if the associated bundle-like metric $g$ makes $\F$ minimal.

Let $\displaystyle{r,s \in \mathbb{N} \cup \{0\}}$ and consider an $(r,s)$-tensor field $T$ on $Q$. Then $T$ is called a basic tensor if the Lie derivative $L_UT$ and the interior product $i_UT$ vanish for any vector field $U$ tangent to $T\F$. In particular, if $T$ is a $(0,0)$ tensor, then $T$ is said to be a basic function. We denote the set of all basic forms by $\Omega^*_B(\F)$ and all basic symmetric $2$-tensors by $S_B^2(\F)$. By the definition of $R^Q$, $\displaystyle{S^Q\in \Omega_B^0(\F)}$ and $\displaystyle{Ric^Q\in S_B^2(\F)}$ are directly calculated\cite{Ton}. However, $\kappa \not\in \Omega^1_B(\F)$ in general.

Neverthless, the $L^2$ orthogonal decomposition of $\Omega^*(M)$ is given by $$\Omega^*(M)=\Omega^*_B(\F)\oplus (\Omega^*_B(\F))^{\perp}$$ by Alvarez-L\'opez\cite{AL}. Here, the $L^2$-inner product of basic differential $r$-forms $\eta$ and $\omega$ is given by \begin{equation}
g_Q(\eta, \omega)_{L^2}=\int_M \eta \wedge *_Q\omega \wedge \x,
\end{equation} where the basic Hodge star operator $*_Q$ is defined by \begin{equation}
*_Q\eta=(-1)^{(q-r)p}*(\alpha\wedge\x)
\end{equation} and $*$ is the Hodge-star operator of the ambient space $M$.

Therefore, we have the decomposition $\kappa=\kappa^{\perp}+\kb$, where $\displaystyle{\kb \in \Omega^1_B(\F)}$ and $\displaystyle{\kappa^{\perp} \in \Omega^1_B(\F)^{\perp}}$. Likewise, we denote $\tb$ by the basic component of $\tau$.

The behaviors of $\kb$ are one of the most important properties in transverse Riemannian geometry. First, the following fact is well-known:
\begin{prop}\cite{AL,Ton}
Let $\F$ be the transversally oriented Riemannian foliation on a compact $(M,\F,g)$. Then, $\kb$ is closed. Moreover, $\kb$ is an exact $1$-form if and only if $\F$ is taut, independent of the choice of $g_Q$. Furthermore, $\F$ is non-taut if and only if any basic $q$-form is exact.
\end{prop}
Thus, we may assume $\displaystyle{\kb=0}$ if $\F$ of a compact $(M,\F,g)$ is taut. Moreover, the formal adjoint of $d_B$(respectively, $\nD_{tr}$) does not coincide with the transverse divergence operator $\dt$ unless $\F$ is taut because of the following theorem:
\begin{thm}\cite{Ton}
Let $M$ of $(M,\F,g)$ be compact and oriented and $\F$ be transversally oriented. Then for a basic $1$-form $\eta$ on $M$,
\begin{equation}
\int_M\dt\eta dV+\int_Mi_{\tb}\eta dV=0.
\end{equation}
\end{thm}

Thus we denote the formal adjoint of $\nD_{tr}$ by $\db$ and $\db$ is defined by
\begin{equation}
\db=i_{\tau_B}+\dt,
\end{equation} which is called the basic divergence of $(M,\F,g)$.
Accordingly, we abuse the notation $\db$ to define the basic codifferntial on $\Omega_B^*(\F)$ of $(M,\F,g)$ by \begin{equation}\db=(-1)^{p(q-r)}*_Q(d_B-\kb\wedge)*_Q,\end{equation} since the basic codifferential $\db$ is the formal adjoint of $d_B$. Therefore the basic Hodge Laplacian on $\Omega_B^*(\F)$ is given by \begin{equation}
\Db=\db d_B+d_B\db
\end{equation} and the basic rough Laplacian is given by
\begin{equation}
\nD^*_{tr}\nD_{tr}=-\left[\sum_{i=1}^q\nD_{\mathbf{e}_i}\nD_{\mathbf{e}_i}\right]+\nD_{\tb}.
\end{equation}

By the difference between $\db$ and $\dt$, the transverse Weitzenb\"ock formula\cite{Jun} on $(M,\F,g)$ also does not coincide with the Weitzenb\"ock formula on a compact non-foliated manifold, unless the given foliation $\F$ is taut.

\begin{thm}(Transverse Weitzenb\"ock formula for basic forms)
Let $\eta$ be a basic $r$-form on a compact $(M,\F,g)$ for $r=1,...,q$. Then the following equality holds:
\begin{equation}
\Db \eta=\nD^*_{tr}\nD_{tr}\eta+Ric^Q\cdot\eta+A_{\tau_B}\eta,
\end{equation} where \begin{equation}Ric^Q\cdot\eta(X_1,...,X_r)=\sum_{j=1}^r \eta(X_1,...,(i_{X_j}Ric^Q)^{\sharp},...,X_r)\end{equation} for basic vector fields $X_1,...,X_r$ on $(M,\F,g)$ and \begin{equation} A_{\tau_B}\eta=L_{\tb}\eta-\nD_{\tb}\eta \end{equation} is a bundle map on $\Omega^r_B(\F)$ extended from the connection $\nD_{X}\tb$ for some $\displaystyle{X \in \Gamma Q}$.

If, moreover, $\eta$ is a basic $1$-form and $\eta^{\sharp}$ is the vector field satisfying $\displaystyle{g_Q(\eta^{\sharp},\cdot)=\eta}$, then we have \begin{equation}Ric^Q\cdot\eta(X)=Ric^Q(X,\eta^{\sharp}).\end{equation}
\end{thm}

Now let $\eta$ be a basic $1$-form on a compact $(M,\F,g)$. Then
\begin{equation}
(\nD_{tr}^*\nD_{tr}\eta)(\eta^{\sharp})=-\dfrac{1}{2}\Db|\eta|^2+|\nD_{tr}\eta|^2
\end{equation} is directly calculated from the definition of the basic rough Laplacian. Hence, Theorem 2.3 induces the following.

\begin{cor}\cite{Jun}(Transverse Bochner formula for basic $1$-forms)
Let $\alpha$ be a basic $1$-form on a compact $(M,\F,g)$.
Then (2.17) in Theorem 2.3 is replaced with
\begin{equation}
-\dfrac{1}{2}\Db|\eta|^2=|\nD_{tr}\eta|^2+Ric^Q(\eta^{\sharp},\eta^{\sharp})+A_{\tau_B}\eta(\eta^{\sharp})-\Db\eta(\eta^{\sharp}).
\end{equation}
\end{cor}

\section{Transverse Ricci solitons}
\subsection{General facts from transverse Ricci flow}
Let $(M,\F,g(t))$ be a $1$-parameter family of bundle-like metrics on a Riemannian manifold with Riemannian foliation satisfying $\displaystyle{g(t)=g_{\F}+g_Q(t)}$. The \textit{transverse Ricci flow} is the transverse geometric flow on $(M,\F,g(t))$ satisfying the following system\cite{LMR}(cf. (1.3)):
\begin{equation}
\begin{dcases}
\dfrac{dg_Q}{dt}=-2Ric^Q\\
\dfrac{dg_{\F}}{dt}=0.
\end{dcases}
\end{equation}

Let us denote $\{\xi_{\alpha}\}_{\alpha=1,...,p}$ by the local orthonormal basis on $T\F$ again. Then the second condition in (3.1) implies that \begin{equation} \dfrac{d\kb}{dt}=0\end{equation} since the Rummler's formula\cite{Rum} on $(M,\F,g)$ is known as the following:
\begin{equation}
d\x(\xi_1,...,\xi_p,\cdot)=\det(g(\xi_{\alpha},\xi_{\beta}))\kappa.
\end{equation} Therefore, we always assume that $\kb$ satisfies (3.2) for later use.

The following fact is the key idea to prove the existence of the solutions for (3.1):

\begin{prop}\cite{BHV}
Let $\widetilde{\nD}$ be the transverse Levi-Civita connection of the transverse metric $\widetilde{g_Q}$ of a compact $(M,\F,\widetilde{g})$. Then there exists a smooth $1$-parameter family of bundle-like metric $\{g(t):t \in [0,T)\}$ with the corresponding $1$-parameter family of transverse Levi-Civita connection $\{\nD^t:t \in [0,T)\}$, solving \begin{equation}
\begin{dcases}
\dfrac{dg_Q}{dt}=-2Ric^Q-L_{X^t}g_Q\\
\dfrac{dg_F}{dt}=0\\
g_Q(0)=\widetilde{g_Q},
\end{dcases}
\end{equation} where $X^t$ is a basic vector field on $(M,\F)$ given by $\displaystyle{X^t=\sum_{i=1}^q\widetilde{\nD}_{\mathbf{e}_i}\mathbf{e}_i-\nD^t_{\mathbf{e}_i}\mathbf{e}_i}$ and $\{\mathbf{e}_i\}$ is a moving frame on $Q$ with respect to $g_Q(t)$.
\end{prop}

\begin{rmk}\normalfont(3.4) is equivalent to (3.1)\cite{BHV}, parallel to the non-foliated case of Deturck's work on Ricci flow and its equivalent geometric flow\cite{Det}.
\end{rmk} Hence, the existence problem of the transverse Ricci flow on a compact manifold with a Riemannian foliation is justified.

\begin{prop}\cite{BHV}
The initial value problem of (3.1) on a compact foliated manifold $M$ has a solution.
\end{prop}

Now assume that $\F$ and $\F'$ are transversally oriented taut foliations on compact Riemannian manifolds $M$ and $M'$ with the bundle-like metric $g$ and $g'$, respectively. A diffeomorphism $\displaystyle{\phi: (M,\F,g) \xrightarrow{} (M',\F',g')}$ is called a \textit{foliated diffeomorphism} if it satisfies $\displaystyle{d\phi(T\F) = T\F'}$
(cf. \cite{JJ}.). Hence, we define the self-similar solution of a transverse Ricci flow by the following:

\begin{defn}\normalfont
A \textit{self-similar solution} $g(t)$ of transverse Ricci flow on a compact foliated manifold $(M,\F)$ is a bundle-like metric satisfying:
\begin{equation}
g_Q(t)=\sigma(t)\phi_t^*g_0,
\end{equation} where $g_0$ is the transverse metric of $g(0)$, $\{\phi_t\}$ is a $1$-parameter family of foliated diffeomorphisms on $(M,\F)$ and $\sigma(t)$ is the linear scale parameter depends only on $t$.
\end{defn}

\subsection{The equation of transverse Ricci solitons}

We calculate the first variation formula of (3.5) for later use.

\begin{lem}
Let $g_Q(t)$ be a $1$-parameter family of transverse metric on $(M,\F,g(t))$ and let $\{\phi_t\}$ be a family of foliated diffeomorphisms on $(M,\F,g(t))$ associated with a basic vector field $X$. Assume that a $1$-parameter family of transverse metric $g_Q(t)$ satisfies (3.5). Then
\begin{equation}
\dfrac{dg_Q}{dt}=\sigma'(t)\phi_t^*g_0+\sigma(t)L_X\phi_t^*g_0,
\end{equation} for some scale-parameter $\sigma(t)$.
\end{lem}
\begin{proof}
The proof is just a direct calculation. (cf.\cite{Top} for the Ricci soliton on a non-foliated manifold.)
\end{proof}

Thus, the following theorem is the precise statement of Theorem 1.2:

\begin{thm}
Let $g_Q(t)$ be a self-similar solution of (3.1) on a compact $(M,\F,g(t))$. Then \begin{equation}
Ric^Q(g_Q(t))+\dfrac{1}{2}L_Xg_Q(t)=\lambda(t)\text{ }g_Q(t),
\end{equation} where $X$ is a basic vector field on $M$ and $\displaystyle{\lambda(t)=\dfrac{1}{2}\dfrac{d}{dt}(\text{log}\sigma(t))}$. Conversely, let $g_Q$ be a transverse metric of a compact $(M,\F,g(t))$ satisfying \begin{equation}
Ric^Q(g_Q)+\dfrac{1}{2}L_Xg_Q=\lambda\text{ }g_Q
\end{equation} for some constant $\lambda$ and let $\displaystyle{\left\{\phi_t|0\leq t<\dfrac{1}{2\lambda}\right\}}$ be the $1$-parameter family of foliated diffeomorphism generating $\dfrac{1}{1-2\lambda t}X$. Then $(1-2\lambda t)\phi_t^*g_Q$ is a self-similar solution of (3.1).
\end{thm}
\begin{proof} The proof is based on the idea for the non-foliated Ricci soliton\cite{Fut, Top}.
Let $g_Q(t)$ be a transverse Ricci soliton on $(M,\F,g(t))$ with the initial data $\displaystyle{g_Q(0)=g_0}$. Then we already have
\begin{equation}
-2Ric^Q(g_Q(t))=\sigma'(t)\phi_t^*g_0+L_{X}\sigma(t)\phi_t^*g_0.
\end{equation}
Thus, we have
\begin{equation}
-2Ric^Q(\sigma(t)\phi_t^*g_0)=\dfrac{d}{dt}(\text{log}\sigma(t))\sigma(t)\phi_t^*g_0+L_{X}\sigma(t)\phi_t^*g_0
\end{equation} from (3.9) by direct calculation.
Hence, a self-similar solution of (3.1) satisfies (3.7) as desired.

Conversely, let $g_Q$ be a transverse metric satisfying (3.8) and let $\displaystyle{\sigma(t)=(1-2\lambda t)}$. Then \begin{equation}
\phi_t^*Ric^Q(g_Q)+\dfrac{1}{2\sigma(t)}\phi_t^*L_{X}(\sigma(t)g_Q)-\lambda\phi_t^*g_Q=0.
\end{equation} By (foliated) diffeomorphism and scale invariance of $Ric^Q$ and the assumption on $\{\phi_t\}$, we have
\begin{equation}
Ric^Q(\sigma(t)\phi_t^*g_Q)+\dfrac{1}{2}L_{Y(t)}(\sigma(t)\phi_t^*g_Q)-\dfrac{\lambda}{\sigma(t)}\sigma(t)\phi_t^*g_Q=0,
\end{equation} where $Y(t)=\dfrac{1}{\sigma(t)}X$ is a time-dependent basic vector field on $M$. Therefore, $\displaystyle{\sigma(t)\phi_t^*g_Q}$ satisfies (3.8) for all $t$. In other words, it is a self-similar solution of (3.1) by (3.5).
\end{proof}

\begin{defn}\normalfont
A transverse metric $g_Q$ on a compact $(M,\F,g)$ is called a \textit{transverse Ricci soliton} if it satisfies (3.8).
\end{defn}
Therefore, a self-similar solution of (3.1) and a transverse Ricci soliton on a compact $(M,\F,g)$ is identified by Theorem 3.6 and Definition 3.7.
In particular, a transverse Ricci soliton is shrinking(respectively, steady or expanding) if it satisfies $\displaystyle{\lambda>0}$(or $\displaystyle{\lambda=0}$, $\displaystyle{\lambda<0}$, respectively). Also, if the basic vector field $X$ of (3.8) satisfies $\displaystyle{X=\nD f}$ for some basic function $f$, we call $g_Q$ a gradient transverse Ricci soliton. Here, such $f$ is called the \textit{potential function} of (3.13). Note that the equation of a gradient transverse Ricci soliton is also defined by
\begin{equation}
Ric^Q+\nD_{tr}df=\lambda\text{ }g_Q,
\end{equation} since $\displaystyle{L_{\nD f}g_Q=2\nD_{tr}df}$.

\begin{rmk}\normalfont
The authors do not know whether Proposition 3.3 is proved or disproved on a complete noncompact foliated manifold. Nevertheless, we define a transverse Ricci soliton $g_Q$ on a complete noncompact foliated manifold by a transverse metric $g_Q$ satisfying (3.8).
\end{rmk}

\subsection{Transverse Ricci solitons and transverse curvature functionals}
In the spirit of \cite{Per}, a gradient Ricci soliton on a compact non-foliated manifold is regarded as a critical point of $\W$(respectively, $\FF$) functional. However, a gradient transverse Ricci soliton on a compact $(M,\F,g)$ may not be regarded as a critical point of the transverse analogy $\W^Q$(respectively, $\FF^Q$) of $\W$(respectively, $\FF$) since the difference between $\db$ and $\dt$ makes the first variation of $\W^Q$ and $\W$ (respectively, $\FF^Q$ and $\FF$) subtly different(cf. Remark 3.15). Therefore, we review Lin's work on $\FF^Q$\cite{Lin} to define $\W^Q$ from $\W$ and $\FF^Q$.

\begin{defn}\normalfont\cite{Lin} Let $f$ be a basic function on a compact $(M,\F,g(t))$ and let $dV$ be a volume form with respect to $g(t)$. Then the $\FF^Q$ functional is defined as:
\begin{equation}
\FF^Q(g_Q,f)=\int_M(S^Q+|\nD f+\tau_B|^2)e^{-f}dV.
\end{equation} In particular, we define $$\lambda^Q=\text{inf}_{f\in\Omega_B^0(\F)}\FF^Q$$ for $f$ satisfying $\displaystyle{\int_Me^{-f}dV=1}$.
\end{defn}
If $\F$ on the given compact $(M,\F,g)$ is taut, then $\FF^Q$ is identified with $\FF$. Thus, $\FF^Q$ is not scale-invariant but foliated diffeomorphism-invariant, parallel to $\FF$\cite{Lin}. Also, the minimizing basic function $f$ of $\FF^Q$ on a compact $(M,\F,g)$ exists since the existence of $\lambda^Q$ is observed in \cite{Lin}.

Similar to the non-foliated case, we also have the normalized version of $\lambda^Q$:

\begin{defn}\normalfont\cite{Lin}
A normalized functional $\overline{\lambda^Q}$ of $\lambda^Q$ on a compact $(M,\F,g)$ is defined as:
\begin{equation}
\overline{\lambda^Q}=(\V(g))^{\frac{2}{q}}\lambda^Q,
\end{equation} where $\V(g)$ is the volume of $(M,\F,g)$.

\end{defn}

To recall the first variation of $\FF^Q$(respectively, $\lambda^Q$ and $\overline{\lambda^Q}$) on $(M,\F,g(t))$, we review the following calculations.

\begin{lem}\cite{Lin}
On $(M,\F,g(t))$, the first variation of the transverse scalar curvature $S^Q(t)$ is given by:
\begin{equation}
\dfrac{dS^Q}{dt}=g_Q\left(-\dfrac{dg_Q}{dt}, Ric^Q\right)+\dt\dt \dfrac{dg_Q}{dt}+\dt d_Btr_Q(\dfrac{dg_Q}{dt}),
\end{equation} where $g_Q$ on symmetric 2-tensor is the inner product on $S_B^2(\F)$ given by \begin{equation}
g_Q(v,w)=\dfrac{1}{4}\sum_{i,j=1}^qv(\mathbf{e}_i, \mathbf{e}_j)w(\mathbf{e}_i,\mathbf{e}_j).
\end{equation}
\end{lem}

\begin{lem}\cite{Lin}
The following equalities are obtained on a compact $(M,\F,g(t))$:
\begin{equation}\int_M \dt\dt \dfrac{dg_Q}{dt} d\mu=\int_M[-g_Q\left(\dfrac{dg_Q}{dt}, \nD_{tr}(df+\kb)\right)+\dfrac{dg_Q}{dt}(\nD f+\tb, \nD f+\tb)]d\mu,\end{equation}
\begin{equation}\int_M \Db tr_Q\left(\dfrac{dg_Q}{dt}\right)d\mu=-\int_Mtr_Q(\dfrac{dg_Q}{dt})(\Db f+|\nD f|^2)d\mu,\end{equation}
\begin{equation}\int_Mg_Q\left(\nD f, \nD \dfrac{df}{dt}\right)d\mu=\int_M\dfrac{df}{dt} (\Db f+|\nD f|^2)d\mu,\end{equation}
\begin{equation}\dfrac{d}{dt}|\nD f+\tau_B|^2=2g_Q\left(\nD\dfrac{df}{dt}, \nD f+\tau_B\right) -\dfrac{dg_Q}{dt}(\nD f+\tau_B,\nD f+\tau_B),\end{equation} where $\displaystyle{d\mu=e^{-f}dV}$.
\end{lem}

\begin{proof}
Since \begin{equation}\nD_{tr}de^{-f}=(-\nD_{tr} df+\nD f\otimes \nD f)e^{-f}\end{equation} implies \begin{equation} \Db(e^{-f})=-(\Db f+|\nD f|^2)e^{-f},\end{equation} (3.19) and (3.20) are straightforward. We refer to Section 2.4 of \cite{CLN} for the detailed calculation. Also, (3.18) and (3.21) are obtained by the same calculation introduced in the same subsection of \cite{CLN}.
\end{proof}

\begin{rmk}\normalfont
Since we have assumed $\kb$ is a closed $1$-form, \begin{equation}
\nD\kb(X,Y)=g_Q(\nD_X \tb, Y)=g_Q(\nD_Y \tb, X)=\nD\kb(Y,X)
\end{equation} is obtained. i.e. $\nD_{tr}\kb$ is a basic symmetric 2-tensor.
\end{rmk}

Therefore, we have the first variation of $\FF^Q$ applying the previous two lemmas:
\begin{thm}\cite{Lin}
On a compact $(M,\F,g(t))$, we have
\begin{equation}
\begin{split}
\left. \dfrac{d}{dt}\right\vert_{t=0}\FF^Q(g_Q,f)&=\int_M g_Q\left(-\dfrac{dg_Q}{dt}, Ric^Q+\nD_{tr}(df+\kb)\right)d\mu\\&+\int_M(-2\Db f-|\nD f|^2+S^Q+|\kb|^2-2\db\kappa_B)\dfrac{d(d\mu)}{dt},
\end{split}
\end{equation}where $\displaystyle{d\mu=e^{-f}dV}$ and $\displaystyle{\dfrac{d(d\mu)}{dt}=\left( \dfrac{1}{2}tr_Q\dfrac{dg_Q}{dt}-\dfrac{df}{dt}\right)e^{-f}dV}$.
\end{thm}
\begin{proof}
Applying Lemma 3.11 and Lemma 3.12, the first variation is straightforward.
\end{proof}

\begin{rmk}\normalfont
Let us denote $\displaystyle{d\mu=e^{-f}dV}$ again.
If we use $\FF$ instead of $\FF^Q$ on an arbitrary compact $(M,\F,g)$, then (3.25) is not calculated since the following calculation is obtained by (3.18) and (3.21): \begin{equation}
\begin{split} &\int_M \left[\dt\dt \dfrac{dg_Q}{dt}+\dfrac{d}{dt}|\nD f|^2+g_Q\left(\dfrac{dg_Q}{dt}, \nD_{tr}(df+\kb)\right)\right]d\mu\\&=\int_M\left[\dfrac{dg_Q}{dt}(\nD f+\tb, \nD f+\tb)+2g_Q\left(\nD\dfrac{df}{dt}, \nD f\right) -\dfrac{dg_Q}{dt}(\nD f,\nD f)\right]d\mu.
\end{split}
\end{equation} Hence, we have to use $\FF^Q$ functional to investigate the corresponding curvature functional of (3.1).
\end{rmk}

To explain the content in Section 4.3, we focus on the behavior of $\lambda^Q$ functional. Note that $\displaystyle{\int_Me^{-f}dV=1}$ implies $\displaystyle{\int_M\left(\dfrac{1}{2}tr_Q\left(\dfrac{dg_Q}{dt}\right)-\dfrac{df}{dt}\right)e^{-f}dV=0}.$ Moreover, we may assume $(-2\Db f-|\nD f|^2+S^Q+|\tb|^2-2\db\kappa_B)$ is a constant on $\lambda^Q$\cite{Lin}. Therefore, \begin{equation}
\int_M\left(\dfrac{1}{2}tr_Q\left(\dfrac{dg_Q}{dt}\right)-\dfrac{df}{dt}\right)(-2\Db f-|d f|^2+S^Q+|\tb|^2-2\db\kappa_B)e^{-f}dV= 0
\end{equation} is obtained from the second term of (3.25). Moreover, we have the following property of $\lambda^Q$.

\begin{lem}\cite{Lin}
On a compact $(M,\F,g)$, $\lambda^Q$ is nondecreasing along (3.1).
\end{lem}
\begin{proof}
By the definition of $\lambda^Q$, we have \begin{equation}
\dfrac{d\lambda^Q}{dt}=-\int_Mg_Q\left(Ric^Q+\nD_{tr}(df+\kb), \dfrac{dg_Q}{dt}\right)e^{-f}dV,
\end{equation} where $f$ is a minimizing function of $\FF^Q$.

Moreover, we may replace $\dfrac{dg_Q}{dt}$ in (3.4) with $\displaystyle{\dfrac{dg_Q}{dt}=-2(Ric^Q+\nD_{tr}(df+\kb))}$ by Remark 3.2, we obtain \begin{equation}
\dfrac{d\lambda^Q}{dt}=\int_M2|Ric^Q+\nD_{tr}(df+\kb)|^2e^{-f}dV
\end{equation} from (3.25). Hence, the lemma is proved by the positive definiteness of $\dfrac{d\lambda}{dt}$.
\end{proof}

In contrast, the first variation of $\overline{\lambda^Q}$ on an arbitrary compact $(M,\F,g)$ may not be nonnegative
\cite{Lin}, because of the behavior of the volume coefficient. Therefore, we review Lin's work on the nondecreasing property of $\overline{\lambda^Q}$ functional on a compact Riemannian manifold with a taut foliation for later use.

\begin{lem}\cite{Lin}
Let $\F$ be a taut foliation on a compact $(M,\F,g)$ and assume $\displaystyle{\lambda^Q \leq 0}$. Then $\overline{\lambda^Q}$ is nondecreasing along (3.1) on $(M,\F,g)$.
\end{lem}
\begin{proof} By direct calculation, we have
\begin{equation}
\dfrac{d\overline{\lambda^Q}}{dt}= 2(\V(g))^{\frac{2}{q}}(-\dfrac{\lambda^Q}{q\V(g)}\int_M(S^Q-\db df)dV+\int_M|Ric^Q+\nD_{tr}(df)|^2e^{-f}dV)
\end{equation} and
\begin{equation}
-\dfrac{\lambda^Q}{q\V(g)}\int_M S^Q dV \geq -\dfrac{1}{q}(\lambda^Q)^2.
\end{equation} Therefore, \begin{equation}
\dfrac{d\overline{\lambda^Q}}{dt} \geq 2(\V(g))^{\frac{2}{q}}(-\dfrac{1}{q}(\lambda^Q)^2+\int_M|Ric^Q+\nD_{tr}(df)|^2e^{-f}dV)
\end{equation}
is derived. Hence, the result is computed by \begin{equation}
\dfrac{d\overline{\lambda^Q}}{dt}\geq 2(\V(g))^{\frac{2}{q}} \int_M|Ric^Q+\nD_{tr}df-\dfrac{1}{q}(S-\Db f)g_Q|^2e^{-f}dV \geq 0,
\end{equation} which is the same as in \cite{Per}.
\end{proof}

Now we introduce $\W^Q$ functional, modified from $\FF^Q$. The definition of $\W^Q$ is an analogy of the $\W$ functional\cite{Per}.

\begin{defn}\normalfont
Let $f$ be a basic function on a compact $(M,\F,g)$ and let $\sigma(t)$ be the scale parameter depends only on $t$. Then the $\W^Q$ functional is defined by
\begin{equation}
\W^Q(g_Q,f,\sigma)=\int_M[\sigma(S^Q+|\nD f+\tau_B|^2)+f-q](4\pi\sigma)^{-\frac{q}{2}}e^{-f}dV.
\end{equation} In particular, we define $$\mu^Q=\text{inf}_{f\in\Omega_B^0(\F)}\W^Q$$ for $f$ satisfying $\displaystyle{\int_M(4\pi\sigma)^{-\frac{q}{2}}e^{-f}dV=1}$.
\end{defn}

The following is the first variation formula of $\W^Q$ functional:
\begin{thm} On a compact $(M,\F,g)$, we have
\begin{equation}
\begin{split}
\left. \dfrac{d}{dt}\right\vert_{t=0}\W^Q &=\int_M g_Q\left(-\sigma \dfrac{dg_Q}{dt}+\dfrac{d\sigma}{dt} g_Q, Ric^Q+\nD_{tr} (df+\kb)-\dfrac{1}{2\sigma}g_Q\right)d\mu \\
&+\int_M\left[\sigma(S^Q-2\Db f-|\nD f|^2+|\kb|^2-2\db\kb)+(f-q-1)\right]\dfrac{d(d\mu)}{dt},
\end{split}
\end{equation}
where $d\mu=(4\pi\sigma)^{-\tfrac{q}{2}}e^{-f}dV$ and $\displaystyle{\dfrac{d(d\mu)}{dt}=\left( \dfrac{1}{2}\dfrac{dg_Q}{dt}-\dfrac{df}{dt}-\dfrac{q}{2\sigma}\dfrac{d\sigma}{dt}\right)e^{-f}dV}$.
\end{thm}
\begin{proof} Although the computation of the variational formula is parallel to the method in \cite{CLN}, we prove the above theorem since we have to calculate the additional twisted terms $\kb$, $\db\kb$, and $\nD_{tr}\kb$ in $\W^Q$.

Therefore, Theorem 3.14 implies the first variation formula of $\W^Q$ by direct calculation.
\begin{equation}
\begin{split}
\left. \dfrac{d}{dt}\right\vert_{t=0}\W^Q&=\int_M\dfrac{d\sigma}{dt}(S^Q+|\nD f+\tau_B|^2)d\mu+\int_M\left(\dfrac{1}{2}tr_Q(\dfrac{dg_Q}{dt})-\dfrac{q}{2\sigma}\dfrac{d\sigma}{dt}\right)d\mu\\
&+\int_M \left[\sigma(S^Q-2\Db f-|\nD f|^2+|\kb|^2-2\db\kappa_B)+(f-q-1)\right]\dfrac{d(d\mu)}{dt}\\
&-\sigma\int_Mg_Q\left(\dfrac{dg_Q}{dt}, Ric^Q+\nD_{tr}(df+\kappa_B)\right)d\mu
\end{split}
\end{equation}

Also,\begin{equation}
\dfrac{1}{2}tr_Q(\dfrac{dg_Q}{dt})-\dfrac{q}{2\sigma}\dfrac{d\sigma}{dt}=g_Q\left(\dfrac{1}{2\sigma}g_Q,\sigma \dfrac{dg_Q}{dt}-\dfrac{d\sigma}{dt} g_Q\right)\end{equation} is obtained by definition and \begin{equation}\int_M(S^Q+|\nD f+\tau_B|^2)d\mu=\int_M(S^Q-\dt(df+\kappa_B))d\mu\end{equation} is obtained by integrating both sides of (3.23).

As $\displaystyle{S^Q-\dt(df+\kb)}$ is the transverse trace of $\displaystyle{Ric^Q+\nD_{tr}(df+\kb)}$, we obtain the first variation formula (3.35) from (3.36).
\end{proof}

For later application in Section 4.2, we now focus on the behaviors and the existence of $\mu^Q$. Contrary to the existence of $\lambda^Q$\cite{Lin}, the existence problem of $\mu^Q$ is not straightforward because the problem is related to the behavior of $\widetilde{\W^Q}$ is the following:

\begin{equation}
\widetilde{\W^Q}=\int_M (S^Qu^2+|\kb|^2u^2+4|\nD u|^2-2\db\kappa_Bu^2-\dfrac{1}{\sigma}u^2\text{log}u^2)dV,
\end{equation} where $u=e^{-\frac{f}{2}}$ for a basic function $f$. However, the minimizer problem of the functional $\displaystyle{\int_M u\text{log}u dV}$ is solven\cite{Rot} and this result gives the existence of minimizing functions of $\widetilde{\W^Q}$. Therefore, we have the existence problem of the minimizing functions for $\W^Q$.

\begin{thm}(The minimizer problem for $\W^Q$)
There is a basic function $f$ which satisfies $\W^Q(g,f,\sigma)=\mu^Q$ on a compact $(M,\F,g)$.
\end{thm}
\begin{proof}
Since the existence problem of the original $\mu$ is proven\cite{ENM}, we only need to show that the minimizing function is basic because other steps of the proof are the same. Let $\{u_k\}_{k=1,2,...}$ be a sequence of basic functions on $M$ which satisfies $\int_M u_k^2=1$ and $u_k \rightarrow u_0$ as $k \rightarrow \infty$. Then we already have $\int u_0^2=1$ and $u_0 \in C^{\infty}(M)$. Now let $V$ be a leaf-directional vector field on $M$. By assumption on $\{u_k\}$, we have $V(u_k)=0$ for all $k=1,2,...$. Therefore, we eventually get \begin{equation}\lim_{k \rightarrow \infty}V(u_k)=V(u_0).\end{equation} This implies that $u_0$, the minimizing function of $\W^Q$, is a basic function.
\end{proof}

Also, $\mu^Q$ satisfies the following properties.

\begin{lem}
$\mu^Q(g_Q, \sigma)$ is a foliated diffeomorphism-invariant and scale-invariant functional on $(M,\F,g)$.
\end{lem}
\begin{proof} Although the proof is similar to the proof of the non-foliated case\cite{Per}, we show the lemma because we use foliated diffeomorphisms instead of flows on a non-foliated manifold.

By the foliated diffeomorphism-invariance of $S^Q$, it is straightforward that $\W^Q$ is foliated diffeomorphism-invariant(cf. $\FF^Q$ is foliated diffeomorphism-invariant\cite{Lin} implies $\W^Q$ is also a foliated diffeomorphism-invariant). Moreover,
\begin{equation}
\mu^Q(cg_Q, c\sigma)=\int_M [c\sigma(\dfrac{S^Q}{c}+\dfrac{|\nD f+\tb|^2}{c})+(f-q)](4\pi\sigma)^{-\frac{q}{2}}c^{-\frac{q}{2}}c^{\frac{q}{2}}dV
\end{equation} is calculated straightforward, since $\displaystyle{S^Q(cg_Q)=\dfrac{S^Q(g_Q)}{c}}$ is obtained.
Therefore, \begin{equation}\mu^Q(cg_Q, c\sigma)=\mu^Q(g_Q, \sigma)\end{equation} as desired.
\end{proof}

\begin{lem} $\mu^Q(g_Q, \sigma)$ is a nondecreasing negative functional on a compact $(M,\F,g)$ which converges to $0$ as $\sigma \rightarrow 0^+$.
\end{lem}
\begin{proof}
The idea is the same as \cite{KL}. However, the formal proof follows below since the ordinary diffeomorphism used in \cite{KL} is replaced with the foliated diffeomorphisms.

Let $f$ be a minimizing function of $\W^Q$, $T$ be the maximal time of the short-time existence of the solution $g_Q$ to (3.1), and let $\sigma=T-t$ since $\dfrac{d}{dt}\sigma=-1.$ It is clear that $\displaystyle{\W^Q \rightarrow 0 \text{ as } t \rightarrow T^-}$ by defition. For this reason we have $\displaystyle{\mu^Q \rightarrow 0 \text{ as } \sigma \rightarrow 0^+.}$

Note that if we assume $\displaystyle{\int_M(4\pi\sigma)^{-\tfrac{q}{2}}e^{-f}dV=1,}$ we have
\begin{equation}
\left. \dfrac{d}{dt}\right\vert_{t=0}\W^Q =\int_M 2\sigma \left|Ric^Q+\nD_{tr} (df+\kb)-\dfrac{1}{2\sigma}g_Q\right|^2 (4\pi\sigma)^{-\tfrac{q}{2}}e^{-f}dV
\end{equation}
by putting $\displaystyle{\dfrac{dg_Q}{dt}=-2(Ric^Q+\nD_{tr}(df+\kappa_B))}$ and $\displaystyle{\dfrac{d\sigma}{dt}=-1}$. Also, we use the minimizer problem, Cauchy-Schwarz inequality, and the symmetric $2$-tensor-trace inequality(cf. Exercise 1.50 in \cite{CLN}) \begin{equation}|A|^2 \geq \dfrac{1}{q}(tr(A))^2\end{equation} for a basic symmetric 2-tensor $\displaystyle{\dfrac{dg_Q}{dt}-\dfrac{1}{2\sigma}g_Q}$ in (3.43). Then we have \begin{equation}\int_M 2\sigma \left|Ric^Q+\nD_{tr} (df+\kb)-\dfrac{1}{2\sigma}g_Q\right|^2 (4\pi\sigma)^{-\tfrac{q}{2}}e^{-f}dV \geq \dfrac{2\sigma}{q}(\FF^2-\dfrac{q}{2})^2\end{equation} for some basic minimizing function $f$. Thus, for $\sigma >0$, It is clear that $\W^Q$ is nondecreasing as $t$ increases. Therefore, we only need to verify the sign of $\left. \W^Q\right\vert_{\sigma=0}(g_Q,f,\sigma)$.
Hence, we obtain
\begin{equation}\left. \W^Q\right\vert_{\sigma=0}(g_Q,f,\sigma) \leq 0,\end{equation} as desired.
\end{proof}

Thus, the critical points of $\mu^Q$(respectively, $\lambda^Q$) do not coincide with the gradient transverse Ricci soliton(cf. (1.5) or (3.13)) because of the first variation formula (3.45) of $\W^Q$(and (3.29) of $\FF^Q$). Nevertheless, the critical points of $\mu^Q$ and $\lambda^Q$ derived by (3.29) or (3.45) are also transverse Ricci solitons. Hence, the critical points of $\mu^Q$(or $\lambda^Q$) should be introduced as another particular case of the transverse Ricci soliton for later use:

\begin{defn}\normalfont
A transverse metric $g_Q$ on $(M,\F,g)$ is a \textit{twisted gradient transverse Ricci soliton} if it satisfies \begin{equation}
Ric^Q+\nD_{tr}(df+\kb)=\lambda \text{ }g_Q
\end{equation} for some constant $\lambda$ and the given $f$ is called the \textit{twisted potential function}.
\end{defn}

On a compact $(M,\F,g)$, note that (3.13), the gradient transverse Ricci soliton, is identified with (3.47) if $\F$ is taut. Therefore, we should compare (3.13) and (3.47) to determine whether $\F$ is taut or non-taut with assuming the existence of the (twisted) potential function $f$ satisfying either (3.13) or (3.47), since the minimizing problems are not proved in the non-taut case.

\section{Properties of transverse Ricci solitons}
\subsection{General formulae}

Throughout this section, we assume the existence of the \textit{(twisted) potential function} $f$, which is the basic function satisfying the (twisted) gradient transverse Ricci soliton equation (3.13) and (3.47) on an arbitrary compact $(M,\F,g)$ since the only unknown cases in the existence problem are the expanding twisted gradient transverse Ricci soliton case and the gradient transverse Ricci soliton case on a compact $(M,\F,g)$ with non-taut $\F$. 

First, parallel to a Ricci soliton on a compact and non-foliated Riemannian manifold, we have the following identities for a gradient transverse Ricci soliton on a compact $(M,\F,g)$.

\begin{thm}
Let $g_Q$ be a gradient transverse Ricci soliton on a compact $(M,\F,g)$. Then the following equations hold:
\begin{equation}
S^Q-\dt df=q\lambda,
\end{equation}
\begin{equation}
d(S^Q+|\nD f|^2-2\lambda f)=0,
\end{equation}
\begin{equation}
dS^Q=2i_{\nD f}Ric^Q,
\end{equation}
\begin{equation}
\Db S^Q=2|Ric^Q|^2-2\lambda S^Q-dS^Q(\nD f-\tb).
\end{equation}
\end{thm}
\begin{proof}
We refer to \cite{ENM} for the calculations of (4.1), (4.2), and (4.3) since the proofs are the same as the case of non-foliated Ricci solitons. We only need to prove (4.4).
By (4.2), we have \begin{equation}
\Db S^Q+\Db |\nD f|^2=2\lambda \Db f.
\end{equation}

On the other hand, we also have \begin{equation}
-\dfrac{1}{2}\Db |\nD f|^2=|Ric^Q-\lambda g_Q|^2-\dfrac{1}{2}dS^Q(\nD f)-\dfrac{1}{2}\tb(|\nD f|^2)
\end{equation} by (2.22) in Corollary 2.4 and the following formulae:

\begin{equation}
|\nD_{tr}df|^2=|Ric^Q-\lambda g_Q|^2=|Ric^Q|^2-\lambda S^Q-\lambda\dt df
\end{equation}

\begin{equation}
Ric^Q(\nD f, \nD f)=\dfrac{1}{2}dS^Q(\nD f)=\dfrac{1}{2}[ \Db df( \nD f)-L_{\tb}df(\nD f)]
\end{equation}

Therefore, we have \begin{equation}
-\dfrac{1}{2}\Db |\nD f|^2=|Ric^Q|^2-\lambda S^Q-\lambda\dt df-\dfrac{1}{2}dS^Q(\nD f)-\dfrac{1}{2}\tb(|\nD f|^2)
\end{equation} since (4.6) and (4.9) are equivalent.

Putting (4.9) in (4.5), we obtain
\begin{equation}
\Db S^Q-2|Ric^Q|^2+2\lambda S^Q+dS^Q(\nD f)-\tb(2\lambda f)+\tb(|\nD f|^2)=0.
\end{equation}

Thus,\begin{equation}
\Db S^Q=2|Ric^Q|^2-2\lambda S^Q-dS^Q(\nD f)+dS^Q(\tb),
\end{equation} is implied by (4.2).
\end{proof}
We also introduce the following identities for twisted gradient transverse Ricci soliton on a compact $(M,\F,g)$, parallel to Theorem 4.1.

\begin{thm}
Let $g_Q$ be a twisted gradient transverse Ricci soliton on a compact $(M,\F,g)$. Then the following equations hold:
\begin{equation}
S^Q-\dt(df+\kb)=q\lambda,
\end{equation}
\begin{equation}
d(S^Q+|\nD f+\tb|^2-2\lambda f)=2\lambda \kb,
\end{equation}
\begin{equation}
dS^Q=2i_{\nD f+\tb}Ric^Q,
\end{equation}
\begin{equation}
\Db S^Q=-2\lambda S^Q-dS^Q(\nD f) + 2|Ric^Q|^2.
\end{equation}
\end{thm}
\begin{proof}
\begin{enumerate}
\item Since \begin{equation}
S^Q-\dt(df+\kb)=tr_Q(Ric^Q+\nD_{tr}(df+\kb))=q\lambda,
\end{equation} the first equation is easily induced.
\item By (2.6), we have \begin{equation}
dS^Q=2\dt(\nD_{tr}(df+\kb)).
\end{equation}
Also, we obtain \begin{equation}
\dt(\nD_{tr}(df+\kb))=\Db(df+\kb)-L_{\tb}(df+\kb)-i_{(\nD f+\tb)}Ric^Q
\end{equation} from (2.17) in Theorem 2.3. Since \begin{equation}
i_{(\nD f+\tb)}Ric^Q=\lambda(df+\kb)-\dfrac{1}{2}d|\nD f+\tb|^2
\end{equation}
by the definition of twisted gradient transverse Ricci soliton and
\begin{equation}
L_{\tb}(df+\kb)= di_{\tb}(df+\kb)=d(g_Q(\tb, \nD f+\tb))
\end{equation} by Cartan's magic formula, we obtain
\begin{equation}
dS^Q=2d\dt(df+\kb)-2\lambda(df+\kb)+d|\nD f+\tb|^2.
\end{equation}
Hence, \begin{equation}
dS^Q+d|\nD f+\tb|^2-2\lambda df=2\lambda \kb,
\end{equation} as desired.
\item As \begin{equation}
dS^Q=2d\dt(df+\kb)-2i_{\nD f+\tb}Ric^Q
\end{equation} by (4.19) and (4.21), the calculation is straightforward.
\item Take $\db$ on the both sides of (4.13).
Then we obtain
\begin{equation}
\Db S^Q=-\Db|\nD f+\tb|^2+2\lambda \db\kb+2\lambda \Db f
\end{equation} by direct calculation.

On the other hand, we have \begin{equation}
|\nD_{tr}(df+\kb)|^2=|Ric^Q-\lambda g_Q|^2=|Ric^Q|^2-\lambda S^Q-\lambda\dt(df+\kb)
\end{equation}
and
\begin{equation}
dS^Q(\nD f+\tb)=\Db(df+\kb)(\nD f+\tb)-L_{\nD f+\tb}(df+\kb)(\nD f+\tb),
\end{equation} similar to the proof of Theorem 4.1.
Therefore, we have
\begin{equation}
\Db S^Q=2|Ric^Q-\lambda g_Q|^2-dS^Q(\nD f+\tb)-\tb(|\nD f+\tb|^2)+2\lambda \db\kb+2\lambda \Db f
\end{equation} putting (2.22), (4.13), (4.25), and (4.26) in (4.24).
Since \begin{equation}
\dfrac{1}{2}\tb(|\nD f+\tb|^2)=\tb(\lambda\kb-\dfrac{1}{2}S^Q+\lambda f)=\lambda|\tb|^2-\dfrac{1}{2}dS^Q(\tb)+\lambda\tb(f),
\end{equation} we also obtain
\begin{equation}
\Db S^Q=2(|Ric^Q|^2-\dfrac{1}{2}dS^Q(\nD f+\tb)+\dfrac{1}{2}dS^Q(\tb)-\lambda S^Q).
\end{equation} Hence, (4.15) is acquired as desired.
\end{enumerate}
\end{proof}

\begin{rmk}\normalfont
Since an expanding twisted gradient transverse Ricci soliton satisfies \begin{equation}
\kb=\dfrac{1}{2\lambda}d(S^Q+|\nD f+\tb|^2-2\lambda f)
\end{equation} on a compact $(M,\F,g)$, $\F$ is a taut foliation and hence the twisted potential function of a twisted gradient transverse Ricci soliton always exists.
\end{rmk}

\subsection{A shrinking transverse Ricci soliton}
In this subsection, we focus on the properties of shrinking transverse Ricci solitons on a compact foliated manifold. The following rigidity holds on a shrinking transverse Ricci soliton, which is used to prove the main result:
\begin{lem}
A shrinking transverse Ricci soliton on a compact foliated manifold is a twisted gradient transverse Ricci soliton.
\end{lem}
\begin{proof} Although the proof is similar to the approach of Perelman\cite{Per}, we show the theorem since $\phi_t$'s are foliated diffeomorphisms.
Let $\{g_Q(t)\}$ be a self-similar solution of (3.1) on a compact manifold $(M,\F)$ with a Riemannian foliation. That is, $\displaystyle{g_Q(t)=(1-2\lambda t)\phi_t^*g_0}$, where $\lambda$ is a positive constant. Then, the definition of transverse Ricci soliton yields:
\begin{equation}
\mu^Q(g_Q(t_1),\sigma(t_1))=\mu^Q(g_Q(t_2),\sigma(t_2))
\end{equation} for any $t_1$ and $t_2$ in the domain $\left[0,\dfrac{1}{2\lambda}\right)$, since the following observation holds, similar to the calculation from \cite{KL}.

Let $\sigma(t)=\dfrac{t_2-ct_1}{1-c}-t$ be a smooth function on $[t_1,t_2]$, where $\displaystyle{c=\dfrac{1-2\lambda t_2}{1-2\lambda t_1}<1}$. By (4.31), we obtain:
\begin{equation}
\mu^Q(g_Q(t_2), \sigma(t_2))=\mu^Q\left(\dfrac{\sigma(t_2)}{\sigma(t_1)}\phi^*g_Q(t_1),\sigma(t_2)\right)=\mu^Q(\phi^*g_Q(t_1),\sigma(t_1))
\end{equation} according to the definition of shrinking transverse Ricci soliton and scale-invariant property of $\mu^Q$ in Lemma 3.21.

On the other hand, by Lemma 3.22, we obtain:
\begin{equation}
\dfrac{d}{dt}\mu^Q \geq 0.
\end{equation}
Here, if there exists $t' \in [t_1,t_2]$ satisfying $\displaystyle{\dfrac{d}{dt}\mu^Q(t')>0}$, then (4.31) is absurd. Therefore, we have:
\begin{equation}
\dfrac{d}{dt}\mu^Q=0
\end{equation} for all $t \in [t_1,t_2]$.
As
\begin{equation}
\left. \dfrac{d}{dt}\right\vert_{t=0}\mu^Q(g_Q,\sigma)=2\sigma\int_M\left|Ric^Q+\nD_{tr}(df+\kappa_B)-\dfrac{1}{2\sigma}g_Q\right|^2(4\pi\sigma)^{-\tfrac{q}{2}}e^{-f}dV=0
\end{equation}
is obtained by (3.43),
\begin{equation}
Ric^Q+\nD_{tr}(df+\kappa_B)=\dfrac{1}{2\sigma}g_Q
\end{equation} is derived for all $t \in [t_1,t_2]$, as desired.
\end{proof}

For later use, we also prove the following theorem:
\begin{thm}
Let $g_Q$ be a shrinking transverse Ricci soliton on a compact $(M,\F,g)$. Then, $\F$ is taut.
\end{thm}
\begin{proof}
By direct calculation from (4.13) in Theorem 4.2, we have \begin{equation}
\kb=\dfrac{1}{2\lambda}d(S^Q+|\nD f+\tb|^2-2\lambda f)
\end{equation} if $g_Q$ is a shrinking twisted gradient transverse Ricci soliton. Therefore, $\F$ is taut if $g_Q$ is a twisted gradient transverse Ricci soliton. Thus, if the transverse metric $g_Q$ of a compact $(M,\F,g)$ is a shrinking transverse Ricci soliton, then $\F$ is taut by Lemma 4.3.
\end{proof} 

\subsection{A steady transverse Ricci soliton}

The properties of steady transverse Ricci soliton are more subtle than the other cases because the following equation \begin{equation}
2\lambda\kb=d(S^Q+|\nD f+\tb|^2-2\lambda f)
\end{equation} is inconclusive to determine the taut condition when the given twisted gradient transverse Ricci soliton is steady. However, the result of this subsection is quite similar to the shrinking transverse Ricci solitons.

At the beginning of this subsection, we investigate the relation between the steady transverse Ricci soliton and the twisted gradient transverse Ricci soliton.
\begin{lem}
A steady transverse Ricci soliton on a compact foliated manifold is a twisted gradient transverse Ricci soliton.
\end{lem}
\begin{proof}
Let $g_Q(t)$ be a steady transverse Ricci soliton on a compact manifold $M$ with a Riemannian foliation $\F$. Then we have $g_Q(t)=\sigma(t)\phi^*_t g_Q(0)$, where $\{\phi_t\}_{t \in [0,\epsilon)}$ is a family of foliated diffeomorphisms on $(M,\F)$ and $\sigma(t)=1$. Therefore, we obtain
\begin{equation}
\lambda^Q(g_Q(t_1))=\lambda^Q(g_Q(t_2))
\end{equation} for any $t_1$ and $t_2$ in the interval $[0,\epsilon)$ since $\lambda^Q$ is foliated-diffeomorphism invariant\cite{Lin}.
On the other hand, \begin{equation}
\dfrac{d}{dt}\lambda^Q \geq 0
\end{equation} holds by Lemma 3.16. As $\displaystyle{\dfrac{d}{dt}\lambda^Q>0}$ for some $t$ implies \begin{equation}
\lambda^Q(g_Q(t_1))<\lambda^Q(g_Q(t_2)),
\end{equation} we should obtain: \begin{equation}
\left. \dfrac{d}{dt}\right\vert_{t=0}\lambda^Q(g_Q)=\int_M|Ric^Q+\nD_{tr}(df+\kappa_B)|^2e^{-f}dV=0
\end{equation} by (3.29). Therefore, \begin{equation}
Ric^Q+\nD_{tr}(df+\kappa_B)=0
\end{equation} holds for all $t \in [t_1,t_2] \subset [0,\epsilon)$, as desired.
\end{proof}

On the other hand,  we have \begin{equation}
0=d(S^Q+|\nD f+\tb|^2)
\end{equation} putting $\lambda=0$ in (4.38). Moreover, we have the following result induced from (4.44):
\begin{lem}
Let $g_Q$ of a compact $(M,\F,g)$ be a steady transverse Ricci soliton. Then
\begin{equation}
S^Q+|\nD f+\tb|^2=0
\end{equation} holds.
\end{lem}
\begin{proof}
Note that
\begin{equation}
S^Q+|\nD f+\tb|^2=C
\end{equation} for some constant $C$, by (4.44). On the other hand,
\begin{equation}
S^Q-\Db f-\db\kappa_B+g_Q(\tb, \nD f)+ |\tb|^2=0
\end{equation} by taking a transverse trace on a twisted gradient transverse Ricci soliton equation.

Subtracting (4.46) from (4.47), we obtain:

\begin{equation}
-\Db f-|\nD f|^2-\db\kappa_B-g_Q(\nD f, \tb)=-C.
\end{equation}

Integrating both sides of (4.48),
\begin{equation}
-\int_M(\Db f+|\nD f|^2)e^{-f}dV=-C\int_Me^{-f}dV
\end{equation} is derived since we acquire

\begin{equation}
\int_M\db\kappa_Be^{-f}dV=-\int_Mg_Q(\nD f ,\tb) e^{-f}dV.
\end{equation}

On the other hand, we have
\begin{equation}
\int_M\Db e^{-f}=-C\int_Me^{-f}dV,
\end{equation} from (4.49), applying (3.23). Therefore, the constant $\displaystyle{C=0}$ as \begin{equation}
\int_M \Db(e^{-f})dV=0
\end{equation} by Theorem 2.2.
\end{proof}

Thus, we have the main result of this subsection. 
\begin{thm}
If $g_Q$ of a compact $(M,\F,g)$ is a steady transverse Ricci soliton, then $\F$ is taut. Moreover, a transverse metric $g_Q$ of a compact $(M,\F,g)$ is a steady transverse Ricci soliton if and only if it is transversally Ricci-flat.
\end{thm}
\begin{proof}
Since the basic laplacian of $S^Q$ of a twisted gradient transverse Ricci soliton is calculated as  \begin{equation}
\Db S^Q=-2\lambda S^Q-dS^Q(\nD f) + 2|Ric^Q|^2,
\end{equation} we obtain the following: \begin{equation}\Db S^Q=-g_Q( \nD f, \nD S^Q )+2|Ric^Q|^2.\end{equation} Thus, let $x$ be a point of $M$ such that $S^Q(x)$ is the minimum value of the transverse scalar curvature $S^Q$.
Then, \begin{equation}\Db S^Q(x)=2|Ric^Q|^2(x) \geq \dfrac{2}{q}(S^Q)^2(x) \geq 0.\end{equation}

On the other hand, since any basic laplacian of a basic function is identified with the ordinary laplacian of the function\cite{KT}, we have the following: \begin{equation}\Db S^Q=\Delta S^Q,\end{equation} where $\Delta S^Q$ is the ordinary laplacian of $S^Q$. Therefore, \begin{equation} \dfrac{2}{q}(S^Q)^2(x) \leq \Db S^Q(x) \leq 0\end{equation} holds as $S^Q(x)$ is assumed to be the minimum value of $S^Q$. In other words, \begin{equation} \text{min}_{x \in M} S^Q=0\end{equation} by (4.55) and (4.57). Thus, we obtain: \begin{equation} S^Q \geq 0.\end{equation}

In contrast, (4.45) implies that \begin{equation} S^Q \leq 0.\end{equation} Therefore, \begin{equation} S^Q=0,\end{equation} is calculated.
From (4.45) and (4.61), we have \begin{equation}|\nD f+\tb|^2\equiv0.\end{equation} Hence, $\F$ is taut by the positive definiteness of (4.62).

Moreover, a steady twisted gradient transverse Ricci soliton is always transversally Ricci-flat since \begin{equation}\tb+\nD f\equiv0 \end{equation} implies \begin{equation} \nD_{tr}(df+\kb)\equiv0,\end{equation}because of (4.62). As the transverse metric $g_Q$ of a transversally Ricci-flat foliation $\F$ of a compact $(M,\F,g)$ is a steady transverse Ricci soliton, the equivalence holds as desired.
\end{proof}
Note that Theorem 4.8 is one of the key idea to prove Theorem 1.3.

\begin{rmk}\normalfont
Theorem 4.7 also implies that any transversally Ricci-flat foliation on a compact Riemannian manifold is always taut.
\end{rmk}

\subsection{Transverse Ricci solitons on a taut foliation}
Now, we investigate the relation between taut transverse Ricci solitons and the (twisted) gradient transverse Ricci solitons to prove Theorem 1.3. The following theorem is the first statement of Theorem 1.3.

\begin{thm}
Let $\F$ of a compact $(M,\F,g)$ be a taut foliation. Then any transverse Ricci soliton $g_Q$ is a gradient transverse Ricci soliton.
\end{thm}
\begin{proof} First, we know that any shrinking transverse Ricci soliton is a twisted gradient transverse Ricci soliton by Lemma 4.4. As a shrinking twisted gradient transverse Ricci soliton satisfies \begin{equation}
\kb=\dfrac{1}{2\lambda}d(S^Q+|\nD f+\tb|^2-2\lambda f),
\end{equation} any shrinking transverse Ricci soliton is in fact gradient. Also, we showed that any steady transverse Ricci soliton is transversally Ricci-flat by Theorem 4.8. Since a transversally Ricci-flat transverse metric is a gradient transverse Ricci soliton with a constant potential function, the theorem is proved for the steady and shrinking transverse Ricci solitons. Therefore, we only need to prove the behavior of expanding transverse Ricci soliton. Let $g_Q$ be an expanding transverse Ricci soliton on $(M,\F,g)$. Since $\F$ is assumed to be taut, it is clear that the following inequality holds:\begin{equation}
\dfrac{d\overline{\lambda^Q}}{dt}\geq 2(\V(g))^{\frac{2}{q}} \int_M|Ric^Q+\nD_{tr}df-\dfrac{1}{q}(S-\Db f)g_Q|^2e^{-f}dV \geq 0.
\end{equation} Also, $\overline{\lambda^Q}$ is constant along (3.1) by Definition 3.10. In other words, we have $\displaystyle{\dfrac{d\overline{\lambda^Q}}{dt}=0}$ for any $t$. Hence, $g_Q$ is a gradient transverse Ricci soliton as the equality of (3.33) holds if and only if $g_Q$ is a gradient transverse Ricci soliton.
\end{proof}

The following result is the second statement of Theorem 1.3.
\begin{thm}
Let $g_Q$ be a transverse metric of a compact $(M,\F,g)$. If $\F$ is non-taut, then $g_Q$ is an expanding transverse Ricci soliton.
\end{thm}
\begin{proof}
If the given transverse metric is a shrinking transverse Ricci soliton, then $\F$ is taut by Theorem 4.5 since we may assume that any shrinking transverse Ricci soliton is twisted gradient and \begin{equation}
\kb=\dfrac{1}{2\lambda}d(S^Q+|\nD f+\tb|^2-2\lambda f)
\end{equation} holds in this case. Moreover, if the given transverse metric is a steady transverse Ricci soliton, then $\F$ is again taut since \begin{equation} |\nD f+\tb|^2\equiv0\end{equation} is induced by Theorem 4.8. Hence, if $g_Q$ is a transverse Ricci soliton, then $g_Q$ must be expanding else a contradiction occurs. 
\end{proof}

The next theorem is the key idea to prove the last statement of Theorem 1.3.
\begin{thm}
Any Riemannian foliation on a compact manifold endowed with a twisted gradient transverse Ricci soliton is taut.
\end{thm}
\begin{proof}
Let $g_Q$ be a twisted gradient transverse Ricci soliton on a compact $(M,\F,g)$. If $g_Q$ is either a shrinking transverse Ricci soliton or an expanding transverse Ricci soliton, then $\F$ is taut since we have \begin{equation}
\kb=\dfrac{1}{2\lambda}d(S^Q+|\nD f+\tb|^2-2\lambda f)
\end{equation} as $\lambda \neq 0$ is assumed. Therefore, let $g_Q$ be a steady transverse Ricci soliton. Then a steady transverse Ricci soliton makes the associated Riemannian foliation $\F$ taut since we already have calculated \begin{equation} |\nD f+\tb|^2\equiv0\end{equation} in Theorem 4.8.
\end{proof}
By Theorem 4.12, there are no non-taut Riemannian foliation $\F$ of a compact $(M,\F,g)$ associated with a twisted gradient transverse Ricci soliton. Therefore, the following corollary is induced from Theorem 4.12.
\begin{cor}
There are no twisted gradient transverse Ricci soliton on a compact $(M,\F,g)$ with the non-exact mean curvature $1$-form. That is, no nontrivial twisted gradient transverse Ricci soliton exists on a compact $(M,\F,g)$.
\end{cor}
\begin{proof}
Suppose not. Then there are no basic functions $h$ on a compact $(M,\F,g)$ satisfying $\kb=dh$. However, if $g_Q$ is either shrinking or steady, then we have \begin{equation}
\kb=\dfrac{1}{2\lambda}d(S^Q+|\nD f+\tb|^2-2\lambda f)
\end{equation} for some nonzero constant $\lambda$. Therefore, we have to assume that $g_Q$ is a steady twisted gradient transverse Ricci soliton but we always have \begin{equation} |\nD f+\tb|^2\equiv0\end{equation} in this case. Thus, a contradiction occurs.
\end{proof}

\subsection{An expanding transverse Ricci soliton}
Although we proved the main results, it is still unknown when an expanding transverse Ricci soliton is a (twisted) gradient transverse Ricci soliton on a compact $(M,\F,g)$. Nevertheless, we have the following rigidity of expanding gradient transverse Ricci solitons.

\begin{thm}
Let $\F$ be a Riemannian foliation with a basic harmonic mean curvature $\kb$ on a compact Riemannian manifold $M$ with a bundle-like metric $g$. If $g_Q$ is an expanding transverse Ricci soliton with a potential function $f \in \Omega_B^0(\F)$, then $\F$ is transversally Einstein.
\end{thm}
\begin{proof} The idea of the proof is similar to the argument in \cite{CZ}.
Note that \begin{equation} S^Q+|\nD f|^2-2\lambda f=C \end{equation} is computed for some constant $C$ by (4.2) in Theorem 4.1.
Subtracting (4.73) from \begin{equation}
S^Q-\dt df=q\lambda,
\end{equation} we obtain the following: \begin{equation} -\dt df-|\nD f|^2=q\lambda-C-2\lambda f.\end{equation}
Integrating both sides of (4.74) with the weighted volume form $e^{-f}dV$, \begin{equation} \int_M g_Q( \tb, \nD f) e^{-f}dV=\int_M\lambda\left(q-\dfrac{C}{\lambda}-2f\right)e^{-f}dV.\end{equation} is calculated by (4.52).
Here, note that the left-hand side of (4.76) is supposed to be zero because we assume that $\displaystyle{\db\kappa_B}=0$. Thus, we have\begin{equation}\int_M\left(q-\dfrac{C}{\lambda}-2f\right)e^{-f}dV=0\end{equation} since $\displaystyle{\lambda<0}$.

On the other hand, we may calculate the below equation: \begin{equation} 2 f \leq q -\dfrac{C}{\lambda} \end{equation} since \begin{equation} 2\lambda f(x_0)+C=S^Q(x_0)=q\lambda+\Db f \geq q\lambda\end{equation} at a maximum point $x_0$ of $f$ on $M$ and $\lambda<0$.

Therefore, \begin{equation}0 \leq \int_M \left(q-\dfrac{C}{\lambda}-2f\right)e^{-f}dV=0\end{equation} holds by (4.78) and \begin{equation} 2f=q-\dfrac{C}{\lambda}\end{equation} is straightforwardly implied.
\end{proof}

Hence, we may assume that an expanding gradient transverse Ricci soliton on a compact foliated manifold is a transverse metric on a transversally Einstein foliation.

\section{Examples of transverse Ricci solitons}
First, we give a trivial example of a transverse Ricci soliton on a compact product manifold.

\begin{eg}\normalfont(Transverse Ricci soliton on product manifolds)
Let $(M_1,g_1)$ be a Riemannian manifold and $(M_2,g_2)$ be a Ricci soliton. Let us denote the Riemannian connection of $(M_1,g_1)$ and $(M_2,g_2)$ by $D^1$ and $D^2$ respectively. Consider a Riemannian product manifold $\displaystyle{(M_1\times M_2,g_1+g_2)}$ with a foliation $\displaystyle{\F=\{M_1 \times \{p\}:p \in M_2 \}}$. The Levi-Civita connection of $M_1 \times M_2$ satisfies the following equation\cite{GHL}.
\begin{equation}D_XY=D_{X_1}^1Y_1+D_{X_2}^2Y_2,\end{equation}
where $X,Y \in \Gamma(TM)$ such that $X=X_1+X_2$ and $Y=Y_1+Y_2$ for $X_1,Y_1 \in \Gamma(TM_1)$ and $X_2, Y_2 \in \Gamma(TM_2)$, respectively.

We observe two facts from $D^1$ and $D^2$. First, the leaves are totally geodesic. This is because $D_XY=D_X^1Y$ for $X,Y \in \Gamma(TM_1)$, we obtain $\pi(D_XY)=0$ for the projection map $\displaystyle{\pi: T(M_1 \times M_2) \xrightarrow{} TM_2}$. Second, $g_2$ is Riemannian. Pick $X \in \Gamma(TM_1)$. We may assume that $Y, Z \in \Gamma(TM_2)$ by the assumption on the connection. Then $\displaystyle{L_Xg_2(Y, Z)=0}$ is directly calculated. Hence, $(M_1\times M_2,g_1+g_2,\F)$ is a taut transverse Ricci soliton.
\end{eg}

Next, we introduce two examples of non-taut transverse Ricci solitons.
\begin{eg}\normalfont(3-dimensional Carri\`ere torus)
We refer to \cite{Car} for the precise construction of this example. Pick a matrix $A \in SL_2(\mathbb{Z})$ whose trace is strictly greater than 2. Let us denote $\rho, \rho^{-1}$ be eigenvalues of $A$ with the corresponding eigenvectors $\dfrac{\partial}{\partial x}$ and $\dfrac{\partial}{\partial y}$, respectively, and consider the quotient $T^2\times \mathbb{R}/[(\mathbf{x},t)]$ of the manifold $T^2\times \mathbb{R}$, where the equivalence relation is defined by $\displaystyle{[(\mathbf{x},t)]=\{(A^t\mathbf{x},t): \mathbf{x} \in T^2 \text{ and } t\in \mathbb{Z} \}}$. Here, $A^t$ denotes $t$-times matrix multiplication of $A$ itself, not the transpose of $A$.  

In this setting, let $M=T^2 \times \mathbb{R}/[(\mathbf{x},t)]:=T^3_{A}$ and $\F$ be the foliation on $M$, whose distribution $T\F$ is spanned by an eigenvector $\dfrac{\partial}{\partial x}$. Note that $T^3_A$ is a Lie group with the multiplication defined by $\displaystyle{(\mathbf{x}_1,t_1) \cdot (\mathbf{x}_2,t_2)=(\mathbf{x}_1+A^{t_1}\mathbf{x}_2,t_1+t_2).}$

Based on the operation on $T^3_A$, we obtain a bundle-like metric \begin{equation}
g=\rho^{2t}dx^2+\rho^{-2t}dy^2+dt^2,
\end{equation} whose induced transverse metric is defined by the following: $\displaystyle{g_Q=\rho^{-2t}dy^2+dt^2}$. The above metric tensor defines a Carri\`ere torus $(T^3_A,\F,g)$ on the foliated manifold $(T^3_{A},\F)$.

Now let us consider a moving frame on $T_A^3$ $$\{\mathbf{e}_1, \mathbf{e}_2, \mathbf{e}_3\}=\left\{\rho^{-t}\dfrac{\partial}{\partial x}, \rho^{t}\dfrac{\partial}{\partial y}, \dfrac{\partial}{\partial t}\right\}$$ with respect to the given bundle-like metric and the dual coframe $\displaystyle{\{\mathbf{e}_i^{\flat}\}_{i=1,2,3}}$ of the given moving frame on $T^3_A$. Applying Slesar's calculation\cite{Sle}, we compute the following result of Levi-Civita connections on the ambient space.
\begin{equation}
D_{\mathbf{e}_i}\mathbf{e}_j=
\begin{cases}
-\text{log}\rho \text{ }\mathbf{e}_3 \text{ if  } i=j=1\\
\text{log}\rho \text{ }\mathbf{e}_1 \text{ if  } i=1, j=3\\
\text{log}\rho \text{ }\mathbf{e}_3 \text{ if  } i=j=2\\
-\text{log}\rho \text{ }\mathbf{e}_2 \text{ if  } i=2, j=3\\
0 \text{ else.}
\end{cases}
\end{equation}

By (5.3), \begin{equation}
g(R(\mathbf{e}_i, \mathbf{e}_j)\mathbf{e}_i, \mathbf{e}_j)=
\begin{cases}
(\text{log}\rho)^2 \text{ if  } i=1, j=2\\
-(\text{log}\rho)^2 \text{ if  } i=1, j=3\\
-(\text{log}\rho)^2 \text{ if  } i=2, j=3
\end{cases}
\end{equation} is derived.

Accordingly, we acquire \begin{equation}
\nD_{\mathbf{e}_i}\mathbf{e}_j=
\begin{cases}
\text{log}\rho \text{ }\mathbf{e}_3 \text{ if } i=j=2\\
-\text{log}\rho \text{ }\mathbf{e}_2 \text{ if } i=2, j=3\\
0 \text{ else.}
\end{cases}
\end{equation} since $\displaystyle{\pi([\mathbf{e}_1,\mathbf{e}_3])=\pi([\mathbf{e}_1,\mathbf{e}_2])=0}$ implies $\displaystyle{\nD_{\mathbf{e}_1}\equiv0}$ by definition of the transverse Levi-Civita connection. Thus, $\displaystyle{i_{\mathbf{e}_1}R^Q=0}$ and the transverse Riemann curvature (and the transverse Ricci curvature) of $T_A^3$ is computed by the following:
\begin{equation} R^Q=g_Q(R(\mathbf{e}_2,\mathbf{e}_3)\mathbf{e}_2,\mathbf{e}_3)=-(\text{log}\rho)^2 <0.\end{equation} \begin{equation} Ric^Q=-(\text{log}\rho)^2(\mathbf{e}_2^{\flat} \otimes \mathbf{e}_2^{\flat} + \mathbf{e}_3^{\flat} \otimes \mathbf{e}_3^{\flat}).\end{equation} Thus, $T_A^3$ is a compact manifold foliated by a transverse Einstein foliation with negative Ricci curvature. Hence, the given foliated manifold is a transverse Ricci soliton.

Note that $\F$ on $T^3_A$ is not taut. Again using (5.3), we directly have $\displaystyle{\tau_B=-\text{log}\rho \mathbf{e}_3,}$ which does not vanish unless $\text{log}\rho=0$ implied by $\displaystyle{tr(A)=2.}$

By Theorem 1.3, $(T^3_A,\F)$ is not a twisted gradient transverse Ricci soliton since $\F$ is non-taut.
\end{eg}

\addresses
\end{document}